\documentclass[
jmfm,
twoside, 
reqno
]
{amsart}

\usepackage{amsmath,amssymb,amsthm}

\usepackage{enumerate}

\usepackage[colorlinks=true, 
pdfcreator={TeXShop for Mac OS X},
unicode,
pdfstartview=FitH,
pdfpagemode=UseNone,
hyperfootnotes=false,
]{hyperref}

\usepackage{xcolor}
\colorlet{linkcolour}{red!50!black}

\hypersetup{
pdftitle={Time-Periodic Solutions of Burgers Equation}, 
pdfauthor={Magnus Fontes and Olivier Verdier},
colorlinks=true,
linkcolor=linkcolour, citecolor=linkcolour, pagecolor= linkcolour,
}

\newcommand{\demph}[1]{\textbf{\textit {#1}}}

\newcommand{\N}{\mathbb{N}}
\newcommand{\R}{\mathbb{R}}
\newcommand{\Z}{\mathbb{Z}}

\newcommand{\dd}{\,\mathrm d}

\newcommand{\bracket}[2]{\left\langle #1,#2 \right\rangle}
\newcommand{\duality}[2]{\bracket{#1}{#2}}
\newcommand{\bracketD}[2]{\bracket {#1}{#2}_{\mathcal{D}',\mathcal{D}}}

\newcommand{\psL}[2]{\left(#1,#2\right)}
\newcommand{\normL}[1]{\abs{#1}}

\newcommand{\norm}[1]{\left\lVert#1\right\rVert}
\providecommand{\abs}[1]{\left\lvert#1\right\rvert}

\newcommand{\dt}[1]{#1_{\sqrt{t}}}
\newcommand{\dx}[1]{#1_{\sqrt{x}}}
\newcommand{\dts}[1]{#1_{\sqrt{t}*}}
\newcommand{\Dx}[1]{#1_x}
\newcommand{\Dxop}{\operatorname{D_x}}
\newcommand{\Dsop}[1]{\operatorname{D}^{#1}}
\newcommand{\Dssop}[1]{\operatorname{D}^{#1}_{*}}
\newcommand{\Dhalfop}{\Dsop{\frac 12}}

\newcommand{\Hilbop}{\operatorname{\mathcal{H}}}
\newcommand{\Burgop}{\operatorname{T}}
\newcommand{\Linop}{\operatorname{\mathcal{L}}}

\newcommand{\Hilbert}{H}
\newcommand{\DTH}{\mathcal{D}(\T,H)}
\newcommand{\dTH}{\mathcal{D}'(\T,H^*)}

\newcommand{\Sob}[1]{\mathrm{H}^{\left(#1\right)}}
\newcommand{\Sobsob}[2]{\mathrm{H}^{\left(#1\right)\left(#2\right)}}
\renewcommand{\H}{\Sob{\frac 12,1}_0}
\newcommand{\Cont}[1]{\mathcal{C}^{#1}}
\newcommand{\Hs}{\Sob{-\frac 12,-1}}
\newcommand{\Leb}[1]{\mathrm{L}^{#1}}
\newcommand{\Four}[1]{\hat{#1}}
\newcommand{\Fourop}{\operatorname{\mathcal{F}}}
\newcommand{\Hilb}[1]{\widetilde{#1}}

\newcommand{\Q}{Q}
\newcommand{\I}{I}
\newcommand{\T}{\mathbb{T}}
\newcommand{\TS}{\T\times \R}

\newcommand{\TxI}{\T\times\I}
\newcommand{\RxI}{\R\times\I}
\newcommand{\IxI}{(0,1)\times\I}

\newcommand{\puls}{2\pi}
\newcommand{\e}{\mathrm{e}}
\newcommand{\im}{\mathrm{i}}
\newcommand{\expn}[1]{\e^{ \im \puls #1 t}}
\newcommand{\expnm}[1]{\e^{- \im \puls #1 t}}

\newcommand{\sgn}{\operatorname{sgn}}
\newcommand{\Id}{\operatorname{\mathrm{Id}}}

\newcommand{\supp}{\operatorname{supp}}
\newcommand{\ind}{\operatorname{ind}}

\newtheorem{thm}{Theorem}

\newtheorem{cor}{Corollary}[section]

\newtheorem{lma}{Lemma}[section]

\newtheorem{prop}{Proposition}[section]

\theoremstyle{remark}
\newtheorem{rem}{Remark}[section]
\theoremstyle{definition}
\newtheorem{df}{Definition}

\newenvironment{pr}[1][\proofname]{
\small
\begin{proof}[#1]\mbox{}}{\end{proof}
\normalsize}
\begin{document}

\title{Time-Periodic Solutions of the Burgers Equation}

\author{Magnus Fontes}
\email{fontes@maths.lth.se}
\address{Box 118, SE-221 00 Lund, Sweden}

\author{Olivier Verdier}
\email{olivier@maths.lth.se}
\address{Box 118, SE-221 00 Lund, Sweden}
\urladdr{http://www.maths.lth.se/\~{}olivier}


\begin{abstract}
We investigate the time periodic solutions to the viscous Burgers equation $u_t -\mu u_{xx} + uu_x = f$ for irregular forcing terms. We prove that the corresponding Burgers operator is a diffeomorphism between appropriate function spaces.
\end{abstract}


\pagestyle{myheadings}
\thispagestyle{plain}
\markboth{M. FONTES AND O. VERDIER}{TIME-PERIODIC SOLUTIONS OF THE BURGERS EQUATION}

\maketitle

\section{Introduction}

The study of the Burgers equation has a long history starting with the seminal papers by Burgers \cite{Burgers}, Cole \cite{Cole}  and Hopf \cite{Hopf} where the Cole-Hopf transformation was introduced. The Cole-Hopf transformation transforms the homogeneous Burgers equation into the heat equation.

More recently there have been several articles dealing with the forced Burgers equation:
\begin{equation}\label{simple-Burger}
u_t - \nu u_{xx} + u u_x = f
\end{equation}
The vast majority treats the initial value problem in time with homogeneous Dirichlet or periodic space boundary conditions (see for instance \cite{Kreiss}).

Only recently has the question of the time-periodic forced Burgers equation been tackled. In most cases \cite{Moser,WeinanE} the authors are chiefly interested in the inviscid limit (the limit when the viscosity $\nu$ tends to zero). The forcing term is usually chosen to take a particular form, for example a sum of products of white noises in time and smooth functions in space \cite{WeinanE,Sinai}. In \cite{Fokas} the space domain is the half line and the Dirichlet boundary conditions are time periodic and analytic.

The closest related work to ours is that of Jauslin, Kreiss and Moser \cite{Moser} in which the authors show existence and uniqueness of a space and time periodic solution of the Burgers equation for a space and time periodic forcing term which is smooth.

In this article we generalise these results and prove that the Burgers operator coming from the Burgers equation is in fact a diffeomorphism between appropriate time periodic anisotropic Sobolev spaces.

More precisely our main result (\autoref{SmoothnessThm}) shows that given a time periodic forcing term in $\Hs$ 
 we have existence and uniqueness of a time periodic solution in $\H$. Furthermore we are able to prove smooth dependence on the forcing term. 

To prove this result we will use a method similar to that of \cite{Fontes-Saksman} which makes extensive use of anisotropic Sobolev spaces. We will also use a modification of the Cole-Hopf transformation to prove uniqueness of the solution.

We prove those results with homogeneous Dirichlet spatial boundary conditions but the results and the proof can be extended to inhomogeneous boundary conditions as well as spatial periodic boundary conditions.

\section{Fractional Calculus}

In this section we recall some well known facts and fix some general notations.  
\subsection{Fourier development}

We denote the \emph{one dimensional torus} by:
\[\T = \R/\Z\]

Let $\Hilbert$ denote a \emph{complex Hilbert space}, then the space of  smooth \emph{Hilbert space-valued periodic test functions} is denoted by:
\[\DTH  = \Cont{\infty} (\T, \Hilbert) \]
endowed with the usual topology of test functions.
Any test function $\varphi \in \DTH$ can be developped in a Fourier series:
$$\varphi = \sum_{n\in \Z} \varphi_n \cdot \expn{n}$$
where $\varphi_n\in\Hilbert$ is defined by:
$$ \varphi_n = \int_{\T} \varphi(t) \cdot \expnm{n} \dd t$$

We denote by $\dTH$ the space of {continuous linear functionals} on $\DTH$. It is naturally isomorphic to the Hilbert-space valued 1-periodic distributions on $\R$. 
 For any periodic distribution $u \in \dTH$ we then have:
$$u = \sum_{n\in\Z} u_n \cdot \expn{n}$$
where $u_n\in\Hilbert^*$ is defined by
$$ \forall \psi \in \Hilbert,\quad \bracket{u_n}{\psi}_{\Hilbert^*,\Hilbert} := \bracketD{u}{\expnm{n}\cdot \psi}$$

\subsection{Fractional Derivatives}

For any positive real number $s$ we may define the fractional derivative of order $s$ in the following way on $\dTH$ :
$$\Dsop{s} u = \sum_{k\in\Z} (\puls \im k)^s u_k \expn{k} = \sum_{k\in\Z} \abs{\puls \im k}^{s} \e^{\im \sgn(k) s \frac{\pi}{2}}u_k \expn{k}$$
where we have used the principal branch of the logarithm. The sign function is defined as follows:
\[ 
\sgn(k) := \begin{cases} \frac{k}{\abs{k}} & \text{if $k\neq 0$}\\
0 &\text{if $k=0$}\end{cases} 
\] 
For $s=0$ we define $\Dsop{0}=\Id$.
$\Dsop{1}$ coincides with the usual differentiation operator on $\dTH$. The familiar composition property also holds: $\Dsop{s} \circ \Dsop{t} = \Dsop{s+t}$ for any $t,s\geq 0$.

The \emph{adjoint operator} of $\Dsop{s}$ is defined by using the conjugate of the multiplier of $\Dsop{s}$:
$$\Dssop{s} u = \sum_{k\in\Z} \abs{\puls \im k}^{s} \e^{-\im \sgn(k) s \frac{\pi}{2}}u_k \expn{k}$$

$\Dsop{s}$ and $\Dssop{s}$ are adjoints in the sense that for any $u\in\dTH$ and $\varphi\in\DTH$:
$$\duality{\Dsop{s} u}{\varphi} = \duality{u}{\Dssop{s} \varphi}$$
and similarly:
$$\duality{\Dssop{s} u}{\varphi} = \duality{u}{\Dsop{s} \varphi}$$

\subsection{Hilbert Transform}

The \emph{Hilbert transform} $\Hilbop$ is defined using the multiplier $-\im\sgn{k}$. For $u\in\dTH$ let 
$$\Hilbop u = \sum_{k\in\Z}-\im\sgn{k}\,u_k\, \expn{k}$$
 Simple computations then give:
\[ \Dhalfop_* = \Dhalfop\circ\Hilbop = \Hilbop\circ\Dhalfop \]

Notice that if $\Hilbert$ is a function space then $\Hilbop$ maps real functions to real functions.
The following properties will be useful in the sequel:  
\begin{gather}\label{dtshilbeq}
\forall u\in\Sob{\frac 12}(\T,\Hilbert) \quad \psL{\Dhalfop u}{\Dhalfop_* \Hilbop u}_{\Leb{2}(\T,\Hilbert)} = - \norm{\Dhalfop u}^2_{\Leb{2}(\T,\Hilbert)} \\
\forall u \in\Leb{2}(\TxI)\quad \Re(\psL{u}{\Hilbop(u)}_{\Leb{2}(\TxI)}) = 0
\end{gather}
where $\Re$ denotes the real part of the expression.

\subsection{Fractional Sobolev Spaces}

We define fractional Sobolev spaces in the following manner, for any $s\in\R$:
$$\Sob{s}(\T,\Hilbert) = \Bigl\{u \in \dTH;\quad \sum_{k\in\Z}\abs{1+k^2}^{s}\norm{u_k}_{\Hilbert}^2 <\infty\Bigr\}$$

Of course $\Sob{0}(\T,\Hilbert) = \Leb{2}(\T,\Hilbert)$. When $s\geq 0$ then for an $u\in\Leb{2}(\T,\Hilbert)$: $u\in\Sob{s}(\T,\Hilbert) \iff \Dsop{s} u \in \Leb{2}(\T,\Hilbert)$. Moreover $\Sob{s}(\T,\Hilbert)$ is then a Hilbert space with the following scalar product:
$$\psL{u}{v} := \psL{u}{v}_{\Leb{2}(\T,\Hilbert)} + \psL{\Dsop{s} u}{\Dsop{s} v}_{\Leb{2}(\T,\Hilbert)}$$

The following classical result holds: $\bigl(\Sob{s}(\T,\Hilbert)\bigr)^* = \Sob{-s}(\T,\Hilbert^*)$.

\subsection{Anisotropic Fractional Sobolev Spaces}

Let $\I$ be an interval in $\R$ and $s\geq 0$. Let $\Sob{s}(\I)$ denote the usual fractional Sobolev space of real-valued s-times differentiable functions on $\I$.  
$\Sob{s}_0(\I)$ is the closure of $\mathcal{D}(\I)$ in $\Sob{s}(\I)$. In that case we have $\bigl(\Sob{s}_0(\I)\bigr)^* = \Sob{-s}(\I)$.
We will also use the following notations, for $\alpha$, $\beta$ nonnegative real numbers:
$$\Sobsob{\alpha}{\beta}(\TxI) = \Sob{\alpha}(\T, \Sob{\beta}(\I))$$
and
$$\Sob{\alpha,\beta}(\TxI) = \Sobsob{\alpha}{0}(\TxI)\cap \Sobsob{0}{\beta}(\TxI)$$
We also introduce $\Sob{\alpha,\beta}_0(\TxI)$ as the closure of $\mathcal{D}(\TxI)$ in $\Sob{\alpha,\beta}(\TxI$). It is clear that $\Sob{\alpha,\beta}_0(\TxI) = \Sobsob{\alpha}{0}(\TxI)\cap\Leb{2}(\T,\Sob{\beta}_0(\I))$. 
Duals of such spaces are denoted as:
\[\begin{split}\Sob{-\alpha,-\beta}(\TxI) := \Bigl(\Sob{\alpha,\beta}_0(\TxI)\Bigr)^*  &= \Sob{-\alpha}(\T,\Leb{2}(\I)) + \Leb{2}(\T,\Sob{-\beta}(\I)) \\ &=\Sobsob{-\alpha}{0}(\TxI)+\Sobsob{0}{-\beta}(\TxI)\end{split}\]

\subsection{Some embeddings}

If $s_k(\xi)$ is the Fourier transform $s_{k}(\xi) = \Four{u}(k,\xi)$ of a distribution $u$ defined on $\TS$, we have the following H\"o{}lder inequality for any $\theta\in [0,1]$:
\begin{multline*}
\int_{\R}\sum_{k \in  \Z} \abs{k}^{2\alpha(1-\theta)} \abs{\xi}^{2\beta\theta} \abs{s_{k}(\xi)}^2 \dd \xi \leq \\
 \left(\int_\R \sum_{k \in  \Z} \abs{k}^{2\alpha} \abs{s_{k}(\xi)}^2 \dd \xi \right)^{1-\theta}\left(\int_\R \sum_{k \in \Z} \abs{\xi}^{2\beta} \abs{s_{k}(\xi)}^2\dd \xi\right)^{\theta}
 \end{multline*}
From this H\"o{}lder inequality we deduce
\begin{equation*}
\Sob{\alpha,\beta}(\TS) \hookrightarrow \Sob{(1-\theta)\alpha}(\T, \Sob{\theta\beta}(\R))
\end{equation*}
So using an extension operator from $\Sob{\theta\beta}(\I)$ to $\Sob{\theta\beta}(\R)$ one can prove the corresponding inclusion:
\begin{equation}\label{interpolationInclusions}
\Sob{\alpha,\beta}(\TxI) \hookrightarrow \Sobsob{(1-\theta)\alpha}{\theta\beta}(\TxI)\end{equation}

For $\alpha = 1/2$ and $\beta = 1$ and $\theta = \frac 13$ we get:
$$\H(\TxI) \subset \Sob{\frac 1 2, 1}(\TxI) \subset \Sobsob{1/3}{1/3}(\TxI)$$
Then the vectorial Sobolev inequalities yield:
\begin{equation}\label{HincludeLfour}  \H(\TxI) \subset \Sobsob{1/3}{1/3} (\TxI)\hookrightarrow \Leb{4}(\T, \Sob{\frac 1 3}(\I)) \hookrightarrow \Leb{4}(\T, \Leb{4}(\I)) = \Leb{4}(\TxI) \end{equation}
Here the injection $\Sobsob{1/3}{1/3}(\TxI)\hookrightarrow\Leb{4}(\T,\Sob{1/3})$ is compact and thus the injection $\H(\TxI)\hookrightarrow\Leb{4}(\TxI)$ is \emph{compact}.

\section{The Burgers Equation}

\subsection{Preliminary Scalings}

For a period $T>0$, a length $L>0$, a non zero constant viscosity $\nu$ and a \emph{time-periodic forcing term} $g$, the Burgers equation is formally defined on $\R/T\Z\times (0,L)$ by:
$$u_t + u  u_x - \nu u_{xx} = g$$

\newcommand{\scal}[1]{\bar{#1}}

For $\scal{t}\in \T$, $\scal{x}\in  (0,1)$ we define:
$$\scal{u}(\scal{t},\scal{x}) := \frac{T}{L}u(\scal{t} T, \scal{x} L)$$

Then $\scal{u}$ is solution of
$$\scal{u}_t +  \scal{u} \scal{u}_x - \mu  \scal{u}_{xx} = f$$
where
\begin{gather*}
f(\scal{t},\scal{x}) = \frac{T^2}{L}g(\scal{t}T,\scal{x}L)\\
\mu  = \frac{\nu T}{L^2}
\end{gather*}

$\dfrac{1}{\mu}$ is often called the \emph{Reynolds number}. The scalings allow us to restrict the study of the Burgers equation to the \emph{normalised domain} $\R/\Z \times (0,1)$.

\subsection{Notations}

In the sequel we will use the following notation:

\begin{align*}
 \I &:= (0,1)\\
 \Q& := \T \times \I\\
 \psL{u}{v} &:= \int_{\Q} u\cdot v\dd t\dd x\\
 \normL{u} &:= \sqrt{\psL{u}{u}}\\
 \H &:= \H(\Q)\\
 \Hs &:= (\H)^*\\
 \Leb{p} &:= \Leb{p}(\Q)
 \end{align*}
 and for $u\in\H$:
 \begin{align*}
 \norm{u} &:= \norm{u}_{\H}\\
 \Hilb{u} &:= \Hilbop u\\
 \dt{u} &:= \Dhalfop u \in\Leb{2}\\
 \dts{u} &:= \Dhalfop_* u \in\Leb{2}\\
 u_x &:= \frac{\partial u}{\partial x} \in \Leb{2}\\
 \text{For $f\in\Hs$}\quad \duality{f}{u} &:= \duality{f}{u}_{\Hs,\H}
 \end{align*}

\subsection{Functional Setting}

By possibly changing the direction of time  we may always assume that $\mu$ is a positive real number. We split the Burgers equation in a linear and a non-linear part by means of the two following operators:

\begin{df}\label{Linopdef}
We define  $\Linop: \H \rightarrow \Hs$ as:
$$\forall v \in \H \quad \duality{\Linop u}{v} := \psL{\dt u}{\dts v} + \mu \psL{u_x}{v_x}$$
\end{df}

\begin{df} The (nonlinear) operator $S : \H \rightarrow \Hs$ is defined as:
$$\forall v \in \H\quad \duality{S(u)}{v} := -\frac 12 \psL{{u^2}}{\Dx{v}}$$
\end{df}
This operator is well-defined since $\H \subset \Leb{4}$ (cf. \eqref{HincludeLfour}).

\begin{df}The \demph{Burgers operator} $\Burgop : \H \longrightarrow \Hs$ is defined by:
\[
\Burgop = \Linop + S
\]
\end{df}

Given $f\in \Hs$ the Burgers equation becomes:
\begin{equation}
\Burgop(u) = f\label{Burgeq}
\end{equation}

\subsection{Main Result}

Here is the main result of this article:

\begin{thm}\label{SmoothnessThm}
The (nonlinear) Burgers operator $\Burgop$ is a diffeomorphism from $\H$ to $\Hs$, i.e. it is a smooth bijection with smooth inverse.
\end{thm}

The main ingredients in the proof of this result are an a priori estimate and the Cole-Hopf transformation. After giving the a priori estimate we will prove existence and then uniqueness. Before that we make some initial observations.

\subsection{Some elementary properties}

If we denote by $\Burgop'(m)$ the derivative of the operator $\Burgop$ at $m\in\H$ then the following holds for any $u$, $v$ in $\H$:
\begin{equation}\label{square-derivative-observation}
\Burgop(u) - \Burgop(v) = \Burgop'\left(\frac{u+v}{2}\right)\cdot(u-v)
\end{equation}
so $\Burgop$ is injective iff $\Burgop'(m)$ is injective for any $m\in\H$. 

We notice that by the inclusion \eqref{HincludeLfour}: $S(u) = u\cdot u_x \in \Leb{4/3} \hookrightarrow \Hs$ and the last inclusion is the adjoint of the inclusion \eqref{HincludeLfour} and is thus compact as well. Since $u\rightarrow u^2$ is continuous from $\Leb{4}$ to $\Leb{2}$ we deduce that $S$ is a non-linear compact operator (that is to say it is continuous and sends bounded sets of $\H$ to relatively compact sets of $\Hs$). Now as a general fact if $S$ is compact and differentiable then $S'(m)$ is a compact linear operator at any point $m\in\H$. 
We collect these elementary observations in the following lemma:
\begin{lma}\label{Scompact}
The nonlinear operator $S$ is compact and for any $m\in\H$ the derivative  $S'(m)$ at the point $m$ is a linear compact operator.
\end{lma}

\section{An a priori Estimate}

\newcommand{\NTS}{N_{\TS}}
\newcommand{\homog}{\mathcal{C}}
\newcommand{\NQ}{N_Q}

\newcommand{\Rez}{R_0}

We have the following a priori estimate of the solution set: 
\begin{thm}\label{lambdaaprioriestimatethm}
Let $f\in\Hs$. The set 
$$\bigcup_{\lambda\in [0,1]} (\Linop +\lambda S)^{-1}\bigl(\{f\}\bigr)$$
is bounded in $\H$.
\end{thm}
We get the following immediate corollary for the case $\lambda = 1$ since $\Burgop = \Linop + S$:
\begin{cor}\label{aprioriestimatethm}
Let
$f\in\Hs$. The set $\Burgop^{-1}(\{f\})$ is bounded in $\H$.  
\end{cor}

To prove \autoref{lambdaaprioriestimatethm} we will use the same techniques as in \cite{Fontes-Saksman}. In particular, the following lemma:

\begin{lma}\label{DistRepLma}
Given $f\in \Hs$ and $\varepsilon>0$ there exists $g$ and $h$ in $\Leb{2}(\Q)$ such that for any $\varphi\in \H$ we have:
\begin{gather}
\normL{g}\leq \varepsilon\\
\langle f, \varphi \rangle = (g, \dts{\varphi}) - (h, \varphi_x) \qquad \forall \varphi\in \H
\end{gather}
\end{lma}

In other words we have 
$$f = \dt{g} + h_x$$
in the distribution sense, and $g$ can be taken as small as we want in $\Leb{2}(\Q)$.

\begin{pr}[Proof of the Lemma]
This follows directly from the fact that $\Sob{0,-1}$ is a dense subspace of $\Hs$. Indeed given an $\varepsilon>0$ there is a $\varphi\in\Sob{0,-1}$ such that $\norm{f-\varphi}_{\Hs}\leq\varepsilon$. By the Hahn-Banach theorem there exist functions $g$, $h_1$ and $h_2$ in $\Leb{2}$ such that $f-\varphi = \dt{g} + {h_1}_x$, $\varphi = {h_2}_x$ and 
$\normL{g}\leq\norm{f-\varphi}\leq\varepsilon$. We take $h = h_1 + h_2$ and the lemma is proved.
\end{pr}

To prove \autoref{lambdaaprioriestimatethm} we will also need the following Gagliardo-Nirenberg type inequality, for which we give an elementary proof for the convenience of the reader:

\begin{lma}\label{ThmHomog}
There exists a constant $\homog\in\R$ such that for any $u\in \H(\Q)$:
$$\int_\Q \abs{u(t,x)}^4 \dd t\dd x \leq \homog^2 \left(\int_\Q \abs{u}^2 \dd t \dd x + \int_\Q \abs{\dt{u}}^2\dd t\dd x\right) \cdot \left(\int_\Q \abs{u_x}^2\dd t\dd x \right)$$
which implies that:
\begin{equation}\label{usqlequux}
\normL{u^2}\leq \homog \norm{u} \normL{u_x}
\end{equation}
\end{lma}

\newcommand{\Rtwo}{\R^2}

\begin{pr}
\begin{enumerate}
\item
Using the standard Sovolev embedding:
$$\Sob{\frac 1 2}(\Rtwo)\subset \Leb{4}(\Rtwo)$$
we get  by a scaling argument:
\begin{multline}\label{sobolevRtwo}
\forall v\in\Sob{\frac 12}(\Rtwo) \\ \int_{\Rtwo} \abs{v(t,x)}^4 \dd t \dd x \leq C \left(\int_{\Rtwo}\abs{\dt{v}(t,x)}^2\dd t \dd x\right) \cdot \left(\int_{\Rtwo}\abs{\dx{v}(t,x)}^2 \dd t \dd x\right)\end{multline}
\item We use the partial Fourier transform in $x$:
$$\Four{v}(t,\xi) = \int_{\R}v(t,x) \e^{-\im \puls x \xi} \dd x$$
$$ \dx{v}(t,x) = \int_{\R} \Four{v}(t,\xi) \sqrt{\im \puls \xi} \e^{\im\puls \xi} \dd \xi$$
By Plancherel and Cauchy-Schwarz:
\begin{equation*}
\begin{split}
\normL{\dx{v}}^2 	&= \int_{\Rtwo} \abs{\Four{v}(t,\xi)} \cdot \abs{\Four{v}(t,\xi)} \abs{\im \puls \xi} \dd \xi \dd t\\
				&\leq  \sqrt{\int_{\Rtwo}\abs{\Four{v}}^2} \cdot \sqrt{\int_{\Rtwo} \abs{\Four{v}(t,\xi)}^2 (\puls \xi)^2 \dd \xi\dd t}\\
				&= \sqrt{\int_{\R^2}\abs{v}^2 \cdot \int_{\R^2}\abs{v_x}^2}
\end{split}
\end{equation*}
\item From the last inequality together with \eqref{sobolevRtwo}, by extending functions by zero outside $\R\times\I$ we get for any $v\in\H(\RxI)$:
$$\int_{\RxI} \abs{v}^4 \leq \int_{\RxI} \abs{\dt{v}}^2 \sqrt{\int_{\RxI} \abs{v}^2 \int_{\RxI} \abs{v_x}^2}$$
The Poincar\'e{} inequality on $\RxI$:
$$\int_{\RxI} \abs{v}^2 \leq \frac{1}{\pi} \int_{\RxI} \abs{v_x}^2$$
then gives
\begin{equation}\label{GNRI}
\int_{\RxI}\abs{v}^4 \leq \frac{C}{\pi} \int_{\RxI}\abs{\dt{v}}^2 \int_{\RxI}\abs{v_x}^2
\end{equation}
\item Finally, given $u\in\H(\Q)$ we define $\tilde{u}$ on $\RxI$ as the only  1-periodic function in $t$ which is equal to $u$ on $(0,1)\times\I$. Take $\varphi$ in the Schwartz space $\mathcal{S}(\R)$ such that 
$\supp(\Four{\varphi})\subset(-\frac 12, \frac12)$, and $\varphi(0) = 1$. Moreover, given $0<\delta<1$, by means of scalings we may always choose $\varphi$ such that $\varphi([-\frac 12, \frac12])\subset[1-\delta, 1+\delta]$. 

\begin{enumerate}
\item $\tilde{u}\in\mathcal{S}'(\R,\Leb{2}(\I))$ so $\Fourop{\tilde u} = \sum_{k\in\Z} u_k \delta(\tau - k)$ where $u_k \in \Leb{2}(\I)$ are the Fourier coefficients of $u$. By the convolution formula we get the following Fourier expansion for $\varphi \tilde{u}$:
$$\Fourop{(\varphi \tilde{u})} = \sum_{k\in\Z}u_k \Four{\varphi}(\tau - k)$$
\item Thus
\begin{equation*}
\begin{split}
\int_{\RxI} \abs{\dt{(\varphi \tilde{u})}}^2 &= \int_{\R} \norm{\dt{(\varphi \tilde{u})}}_{\Leb{2}(\I)}^2 \\ &=\int_{\R} \abs{\tau} \norm{\sum_k u_k \Four{\varphi}(\tau-k)}_{\Leb{2}(\I)}^2\dd\tau \\
&= \sum_k \norm{u_k}^2_{\Leb{2}(\I)}\int_{k-\frac12}^{k+\frac12}\abs{\tau}\abs{\Four\varphi(\tau-k)}^2\dd\tau\\
\end{split}
\end{equation*}
Now the term on the right hand side can be estimated as follows:
\begin{equation*}
\begin{split}
\int_{k-\frac12}^{k+\frac12}\abs{\tau}\abs{\Four{\varphi}(\tau-k)}^2\dd\tau &=
	\int_{k-\frac12}^{k+\frac12}\underbrace{(\abs{\tau} - \abs{k})}_{\leq  1}\abs{\Four{\varphi}(\tau-k)}^2 + \abs{k}\abs{\Four{\varphi}(\tau-k)}^2 \dd \tau \\
	&\leq (1+\abs{k})\int_{\R}\abs{\Four{\varphi}}^2
\end{split}
\end{equation*}
so we get:
\begin{equation}\label{RIphit} \int_{\RxI} \abs{\dt{(\varphi \tilde{u})}}^2 \leq \int_{\R}\abs{\Four\varphi}^2 \cdot \left(\int_{\TxI}\abs{u}^2 + \int_{\TxI} \abs{\dt{u}}^2 \right) \end{equation}
\item  Furthermore
\begin{equation}\label{RIphifour} \int_{\RxI} \abs{\varphi \tilde{u}}^4 \geq (1-\delta)^4\int_{\TxI} \abs{u}^4 \end{equation}
\item Since $\varphi\in\mathcal{S}(\R)$ there exists an $A>0$ such that $\forall t \in \R\quad\abs{\varphi(t)}\leq A/(1+t)$. Thus finally:
\begin{equation}
\begin{split}
\int_{\RxI}\abs{(\varphi\tilde{u})_x}^2 &= \sum_{k=0}^{\infty} \int_{[k,k+1]\cup[-k-1,-k]} \abs{\varphi(t)}^2 \norm{\tilde{u}_x(t)}_{\Leb{2}(\I)}^2\\
&\leq \sum_{k=0}^{\infty} 2\left(\frac{A}{k+1}\right)^2\int_{\Q}\abs{u_x}^2\\
&= C \int_{\Q}\abs{u_x}^2 \label{RIphix}
\end{split}
\end{equation}
\end{enumerate}
\item By using \eqref{GNRI} with $v = \varphi \tilde{u}$ and combining \eqref{RIphit}, \eqref{RIphifour} and \eqref{RIphix}  we get the desired inequality.
\item By using the Poincar\'{e} inequality once more one gets \eqref{usqlequux}, which concludes the proof of \autoref{ThmHomog}.
\end{enumerate}
\end{pr}

We are now ready for the proof of the a priori estimate.

\begin{pr}[Proof of \autoref{lambdaaprioriestimatethm}]
By definition $\Linop u + \lambda S(u) = f$ means:
\begin{equation}
\forall v\in \H\quad \psL{\dt{u}}{\dts{v}} + \mu\psL{u_x}{v_x} -\frac 12 \lambda\psL{u^2}{v_x} = \langle f, v \rangle
\label{burgersdefvareq}
\end{equation}
\begin{enumerate}
\item
We notice that for smooth $u$:
\begin{align*}
\psL{u^2}{u_x} &= \int_{\Q} u^2 u_x\\
				&= \frac{1}{3}\int_{\Q} (u^3)_x \\
				&= 0
\end{align*}
and then by density and continuity this holds for all $u\in\H$.
\item
With $v=u$ in \eqref{burgersdefvareq} we get:
$$\underbrace{\psL{\dt{u}}{\dts{u}}}_{=0} + \mu \psL{u_x}{u_x} + \frac{1}{2}\lambda\underbrace{\psL{u^2}{u_x}}_{=0} = \langle f, u \rangle$$
which gives:
\begin{align}
\normL{u_x}^2&= \frac{\langle f, u \rangle}{\mu }\label{muuxequalf}\\
			&\leq \frac{\norm{f}\norm{u}}{\mu}\label{gooduxestimate}
\end{align}
\item
With $v=\Hilb{u}$ in \eqref{burgersdefvareq} we get:
$$\psL{\dt{u}}{\dts{\Hilb u}} + \mu \underbrace{\psL{u_x}{\Hilb{u}_x}}_{=0} + \frac{1}{2} \lambda \psL{u^2}{\Hilb{u}_x} = \langle f, \Hilb{u} \rangle$$
Using the identity \eqref{dtshilbeq}, the fact that $\norm{\Hilb{u}} = \norm{u}$ and that $\lambda \leq 1$ we get:
\begin{equation}
\normL{\dt{u}}^2\leq \frac{1}{2}\abs{\psL{u^2}{\Hilb{u}_x}} + \norm{f}\norm{u}
\label{dtuestimate}
\end{equation}
\item
We estimate $\abs{\psL{u^2}{{\Hilb{u}}_x}} $ using the \autoref{ThmHomog}:
\begin{equation}\begin{split}
\abs{\psL{u^2}{{\Hilb{u}}_x}} &\leq  \normL{u^2}\normL{u_x} \\
				&\leq  \homog \norm{{u}} \normL{u_x}^2 \label{rawconvectionestimate}
\end{split}\end{equation}
\item
We use the \autoref{DistRepLma} to write $f = \dt{g} + h_x$ together with \eqref{muuxequalf} we have:
\begin{align}
\mu\normL{u_x}^2 &= \langle f, u \rangle\nonumber\\
			&= \psL{g}{\dts{u}} - \psL{h}{u_x} \nonumber\\
			&\leq  \normL{g}\normL{\dt{u}} + \normL{h}\normL{u_x} \nonumber\\
			&\leq  \normL{g}\normL{\dt{u}} + \normL{h}\sqrt{\frac{\norm{f}{\norm{u}}}{\mu}} \label{specialuxestimate}
\end{align}
\item
Using the estimate \eqref{specialuxestimate} inside \eqref{rawconvectionestimate} and the fact that $\normL{\dt{u}}\leq \norm{u}$ we get:
$$\frac 12\abs{\psL{u^2}{\Hilb{u}_x}}\leq \Rez \left[ \normL{g}\norm{{u}}^2 + \normL{h}\sqrt{\frac{\norm{f}}{\mu}}\norm{u}^{\frac{3}{2}} \right] $$
Where $\Rez$ is defined as:
$$\Rez = \frac{\homog}{2} \cdot \frac{1}{\mu}$$

So if we choose $g$ small enough such that 
$$ \Rez  \normL{g}\leq \frac{1}{2}$$
then using \eqref{dtuestimate} we get:
\begin{equation}\label{dtuRoestimate}
\normL{\dt{u}}^2 \leq  \norm{f}\norm{u} +  \frac 1 2 \norm{{u}}^2 + \Rez\normL{h}\sqrt{\frac{\norm{f}}{\mu}}\norm{u}^{\frac{3}{2}}
\end{equation}
So with the notations:
\begin{align*}
a &= 2\left(1 + \frac{1}{\mu}\right)\norm{f}\\
\intertext{and}
b &= 
\Rez\normL{h}\sqrt{\frac{\norm{f}}{\mu}}
\end{align*}
from \eqref{gooduxestimate} and \eqref{dtuRoestimate} we get
$$\norm{u}^2 \leq  a\norm{u} + 2b \norm{u}^{\frac 3 2}$$
A straightforward computation leads to the bound
$$\norm{u}\leq (b+\sqrt{a+b^2})^2$$ 
Since that estimate does not depend on $\lambda$ the theorem is proved.
\end{enumerate}
\end{pr}

\section{Existence of solutions}

\subsection{Existence and Uniqueness in the Linear Case}

\begin{thm}\label{linearcasethm}
$\Linop$ is a continuous bijection.
\end{thm}

\begin{pr}
\begin{enumerate}
\item
Let us define the following operator on $\H$:
\newcommand{\Pop}{\operatorname{P}}
$$\Pop(u) = \frac{u - \Hilb{u}}{\sqrt 2}$$
$\Pop$ is an isomorphism on $\H$ since the corresponding Fourier multiplier has either module one or $1/\sqrt{2}$.
\item Now
\begin{equation*}
\begin{split}
\duality{\Pop^* \Linop u}{u} &= \duality{\Linop u}{\Pop u}\\
					&= \frac{1}{\sqrt 2}\Biggl(\underbrace{\psL{\dt{u}}{\dts{u}}}_{=0} - \psL{\dt{u}}{\dts{\Hilb u}} + \mu \psL{u_x}{u_x} - \mu \underbrace{\psL{u_x}{\Hilb{u_x}}}_{=0} \Biggr)\\
					&= \frac{1}{\sqrt 2} \left(\normL{\dt u}^2 + \mu  \normL{\Dx{u}}^2\right)\\
					&\geq  \frac{\min \{1,\mu\}}{\sqrt{2}} \norm{u}^2
\end{split}
\end{equation*}
\item
$\Pop^* \circ \Linop$ is therefore a coercive and continuous linear operator from $\H$ to $\Hs$. By the Lax-Milgram theorem (cf. for example \cite{Brezis}) it is invertible. Since $\Pop$ is an isomorphism, so is $\Pop^*$. We conclude that $\Linop$ is an isomorphism.
\end{enumerate}
\end{pr}

\subsection{Existence for the General Case}

\begin{thm}\label{nonlinearsurjectivity}
$\Burgop$ is surjective.
\end{thm}
\begin{pr}
\begin{enumerate}
\item 
Because of \autoref{linearcasethm}, 
the equation $\Burgop (u) = f$ can be rewritten:
$$ \Bigl[\Id + K\Bigr] (u) = \Linop^{-1} f$$
where $$K = \Linop^{-1}\circ S$$
so we only have to show that the application $\Id + K$ is surjective. But since $ \Linop^{-1}$ is continuous and by \autoref{Scompact} $S$ is compact, $K$ is a compact map. 
\item
We choose an open ball $U$ of $\H$ that contains $\bigcup_{\lambda\in [0,1]} (\Linop +\lambda S)^{-1}(f)$. \autoref{lambdaaprioriestimatethm} ascertains that for all $\lambda\in[0,1]$, $\Linop^{-1}f\not\in(\Id+\lambda K)(\partial U)$. So the Leray-Schauder degree of $\Id+K$ on $U$ is equal to the one of $\Id$ which is one:
$$D(\Id+K,U,\Linop^{-1}f) = D(\Id, U, \Linop^{-1}f) = 1$$
 As a result, $\Id+ K$ is surjective and the theorem is proved. (For the Leray-Schauder degree theory, see for instance \cite{Deimling}).
\end{enumerate}
\end{pr}

\newcommand{\HW}{\Sob{1,2}_{N}}
\newcommand{\Hphi}{\Sob{1,2}_{N+}}
\newcommand{\HWq}{\HW/\R}
\newcommand{\Hphiq}{\Hphi/\R_+}

\section{Uniqueness}

\begin{thm}\label{uniquenessThm}
$\Burgop$ is injective.
\end{thm}

The proof is quite long and involved and we shall split it into three propositions.

The first observation we make is that if two functions $u,v\in\H$ satisfy $\Burgop(u)=\Burgop(v)$ then $w:=u-v$ satisfies
\begin{equation}\label{BurgdiffeqRk}
\Burgop(w) = -(vw)_x
\end{equation}
Thus to prove \autoref{uniquenessThm} it is suffices to prove that given any fixed $v\in\H$ the equation \eqref{BurgdiffeqRk} has only the trivial solution $w=0$ in $\H$.

This will be done using the Cole-Hopf transformation. In order to define it we shall need the following anisotropic Sobolev space with Neumann boundary conditions:
$$\HW = \Bigl\{u\in \Sob{1,2}(\Q);\ u_x(t,0) = u_x(t,1) = 0 \quad \forall t \in \T \Bigr\}$$
Notice that by \eqref{interpolationInclusions} we have: 
\begin{equation}\label{HWsubsetCzero}
\Sob{1,2}(\Q) \subset \Sobsob{\frac 23}{\frac 23}(\Q) \subset \Cont{0}\bigl(\overline{\T},\Cont{0}(\overline{\I})\bigr) = \Cont{0}(\overline{\Q})
\end{equation}
We may therefore define the following set:
$$\Hphi = \Bigl\{ u\in\HW;\ u > 0\text{ on }\overline{Q}\Bigr\}$$

We will also need the quotient sets $\HWq$ and $\Hphiq$
where the latter is the quotient with respect to the action of the multiplicative group $(\R_+, \times )$ given by the scalar multiplication (i.e. $\varphi\sim \psi \iff \exists \eta>0\text{ s.t. } \psi=\eta\varphi$).

We now define the following three solution sets, all depending on a fixed function $v\in\H$.

\newcommand{\Solone}{S_1}
\newcommand{\Soltwo}{S_2}
\newcommand{\Solthree}{S_3}

\begin{df} We say that $w\in\Solone$ if $w\in\H$ and 
\begin{equation*}
\Burgop(w) = -(vw)_x
\end{equation*}
\end{df}

\begin{df}
We say that $([W],K)\in\Soltwo$ if $[W]\in\HWq$, $K\in \R$ and 
\begin{equation*}
W_t - \mu W_{xx} + \frac 12 (W_x)^2 = -vW_x + K
\end{equation*}
\end{df}

\begin{df}
We say that $([\varphi],K)\in\Solthree$ if $[\varphi]\in\Hphiq$, $K\in\R$ and
\begin{equation}\label{SPeq}
\varphi_t -\mu\varphi_{xx} + v\varphi_x + K\varphi = 0
\end{equation}
\end{df}

Notice that the definitions above are consistent since the equations used do not depend on the chosen representative.
By the remark above, \autoref{uniquenessThm} will be proved if we can show that the cardinality of $\Solone$ is one. We will do this by first proving that the cardinalities of $\Solone$ and $\Solthree$ are the same (\autoref{cardinalityprop}) and then finally by proving that $\operatorname{card}(\Solthree) = 1$ (\autoref{descriptionSolthreeLma}).

We first prove an auxiliary lemma that will be used to construct a bijection between $\Soltwo$ and $\Solthree$:
\begin{lma}\label{exploglma}
The exponential function  is a bijection from $\HW$ to $\Hphi$. The natural logarithm  is its inverse. These functions can be naturally extended to bijections between $\HWq$ and $\Hphiq$.
\end{lma}
\begin{pr}
\begin{enumerate}
\item
Take $W\in\HW$. By \eqref{HWsubsetCzero} we have $W\in\Cont{0}(\overline{Q})$ and thus $\exp(W)\in\Cont{0}(\overline{Q})$. The lemma thus follows by simple computations and Sobolev injections. Indeed since $W_x\in\H$ and $\H\subset\Leb{4}$ (by \eqref{HincludeLfour}), by considering $(\exp W)_{xx}$ we obtain:
$$(\exp W)_{xx} = \Bigl(\underbrace{W_{xx}}_{\in \Leb{2}}+ (\underbrace{W_x}_{\in \Leb{4}})^2\Bigr)\underbrace{\exp W}_{\in \Cont{0}(\overline{\Q})} \quad \in \Leb{2}$$
and since $(\exp W)_t\in\Leb{2}$ we get $\exp(W)\in \Hphi$. 
\item
The proof goes along the same lines for the logarithm function.
\item The exponential and logarithm functions preserve the group actions used to define the quotient sets $\HWq$ and $\Hphiq$ and can thus be extended to bijections between those sets.
\end{enumerate}
\end{pr}

We are now ready to prove our first proposition:
\begin{prop}\label{cardinalityprop}
The cardinalities of the solution sets $\Solone$, $\Soltwo$ and $\Solthree$ defined above are the same.
\end{prop}
\begin{pr}
\begin{enumerate}
\item We shall explicitely construct two transformations, one from $\Solone$ to $\Soltwo$: $\phi_{21}: \Solone\rightarrow\Soltwo$ and the other from $\Soltwo$ to $\Solone$: $\phi_{12}: \Soltwo\rightarrow\Solone$ that are inverse to each other.

For $w\in\Solone$ we have $w_t = (\mu w_{x} - \frac 12 w^2 - vw)_x$ so $w_t\in\Sob{0,-1}$ and thus $\overline{W} := \int_0^x w(t,y)\dd y \in \Sob{1,2}$. Moreover since $w\in\Solone$ we get
$$\Dxop\Bigl(\overline{W}_t  - \mu \overline{W}_{xx} +\frac 12 (\overline{W}_x)^2 + v \overline{W}_x \Bigr) = 0$$
and thus
$$\overline{W}_t - \mu \overline{W}_{xx} + \frac 12 (\overline{W}_x)^2 + v \overline{W}_x = g$$
for some $g\in\Leb{2}(\T)$.

Now define $h\in\Sob{1}(\T)$ by $g(t) = h'(t) + \int_{\T}g$ and then define $W:= [\overline{W} - h]$ i.e. the projection of $\overline{W} - h \in \HW$ onto $\HWq$. An elementary computation shows that $(W,\int_{\T} g)\in\Soltwo$. We put $\phi_{21}(w) = (W, \int_{\T}g)$.
\item On the other hand, given $([W],K)\in \Soltwo$ with $W\in\HW$, a straightforward computation shows that $W_x\in\Solone$. Moreover $W_x$  obviously does not depend on the chosen representative $W$. We put $\phi_{12}([W],K) = W_x$.

Now $\phi_{12}\circ\phi_{21} = \Id$ since for any $h\in\Leb{2}(\T)$:
$$\Dxop\left(\int_0^x w(t,y)\dd y - h(t)\right) = w$$
Furthermore given $([W],K)\in\Soltwo$ then $\phi_{21}\circ\phi_{12}([W],K) = ([U],\tilde{K})\in\Soltwo$ for some $U\in\HW$ and $\tilde{K}\in\R$. One checks that $[U] = [W - h(t)]$ for an $h\in\Sob{1}(\T)$, which implies that $h'(t) = \tilde{K}-K$ and thus since $h$ is periodic $K = \tilde{K}$ and $h$ is constant, and thus $[W] = [W-h]$. As a result $\phi_{21}\circ\phi_{12} = \Id$.
\item We will again construct transformations between the two sets $\Soltwo$ and $\Solthree$ which are inverse to each other: $\phi_{32}: \Soltwo\rightarrow\Solthree$ and $\phi_{23}: \Solthree\rightarrow\Soltwo$.
Given an element $([W],K)\in\Soltwo$ and one representative $W$ one defines $\varphi=\e^{-\frac{W}{2\mu}}$. Then 
\[
\left\{
\begin{aligned}
{\varphi}_t &= -\frac{{W}_t}{2\mu}\varphi\\
{\varphi}_x &= -\frac{{W}_x}{2\mu}{\varphi}\\
{\varphi}_{xx} &= \left[-\frac{{W}_{xx}}{2\mu}+\left(\frac{{W}_x}{2\mu}\right)^2\right]{\varphi}
\end{aligned}\right.
\]
 So
\[
\begin{split}\varphi_t -\mu\varphi_{xx}+v\varphi_x &= -\frac{1}{2\mu}\left[{W}_t -\mu{W}_{xx} +\frac 12 ({W}_x)^2 + v{W}_x\right]\varphi\\
&= -\frac{K}{2\mu}\varphi
\end{split}
\]
Using that $\varphi =  \exp(-W/2\mu)$ and  \autoref{exploglma} one gets $\varphi\in\Hphi$. We denote $[\varphi]$ the projection of $\varphi\in\Hphi$ onto $\Hphiq$. So $([\varphi],\frac{K}{2\mu})$ is in $\Solthree$. Since $[\varphi]$ does not depend on the chosen representative $W$, the function $\phi_{32}$ which maps $([W],K)$ to $([\varphi],K/2\mu)$ is well defined.

\item We define in the same way $\phi_{23}$ by: $\phi_{23}(([\varphi],K)) = ([\log(\varphi)], 2\mu K)$. It is easy to see that $\phi_{23}$ is well defined, that it maps $\Solthree$ to $\Soltwo$, and that $\phi_{23}$ and $\phi_{32}$ are inverse to each other.
\end{enumerate}
\end{pr}

We sum up the last result in the following corollary.
\begin{cor}The Cole-Hopf transformation $\Phi$ defined by $\Phi := \phi_{32}\circ\phi_{21}$ is a bijection from $\Solone$ to $\Solthree$.
\end{cor}

\newcommand{\phim}{\chi}

We will now set out to prove that $\Solthree = \bigl\{([1],0)\bigr\}$. The first step consists in proving that if $(\varphi,K)\in\Solthree$ then $K=0$. We will need a preliminary lemma which proves the positivity of the evolution operator associated with the equation $\psi_t-\mu\psi_{xx}+v\psi_x = 0$. 

\begin{lma}\label{Upositive-lma}
Given $v\in\H$, for any $\psi\in\HW(\IxI)$ such that
\begin{equation}\label{eq-upositive}
\psi_t - \mu \psi_{xx} + v\psi_x = 0
\end{equation}
the traces $\psi(0) := (x\mapsto\psi(0,x))$ and $\psi(1):=(x\mapsto\psi(1,x))$ are well defined in $\Cont{0}(\overline{\I})$ and the following holds:
$$\psi(0)\geq 0 \implies \psi(1)\geq 0$$
\end{lma}
\begin{pr}
\begin{enumerate}
\item The traces are well defined by the same argument as in \eqref{HWsubsetCzero}:  
$$\Sob{1,2}(\IxI)\subset\Sobsob{\frac 23}{\frac 23}(\IxI)\subset \Cont{0}(\overline{\IxI})$$
\item 
\eqref{eq-upositive} implies that
$$\int_{\I} \psi_t \chi \dd x + \mu\int_{\I}\psi_x\chi_x \dd x + \int_{\I} v \psi_x \chi \dd x = 0$$
for any test function $\chi\in\Sob{1,2}(\IxI)$. Using $\chi = \psi^{-} := \max (-\psi(t,x),0)$ we get:
\begin{equation}\label{pairingphiminus}
\int_{\I} \chi_t \chi\dd x + \mu\int_{\I}\chi_x^2 \dd x + \int_{\I} v \chi_x \chi \dd x = 0
\end{equation}
\item We define
\begin{gather*}
g(t) =  \int_{\I} \chi(t,x)^2 \dd x\\
h(t) =  \int_{\I} \chi_x(t,x)^2 \dd x
\end{gather*}
So we can rewrite \eqref{pairingphiminus} as:
\begin{equation}\label{pairingphiminus-gh}
\frac 12 g'(t) + \mu h(t) - \frac{1}{2}\int_{\I} v_x \chi^2\dd x=0
\end{equation}
\item We estimate the third term of the last equation:
\begin{equation}\label{vxchisquare-estimate}
\int_{\I} v_x \chi^2 = \int_{\I} v_x \left(\chi^2 - g(t)\right)\dd x + \int_{\I} v_x g(t) \dd x
\end{equation}
We estimate the two last terms in the following way:
\begin{enumerate}
\item
$$\abs{\int_{\I} v_x (\chi^2 - g(t))\dd x} \leq  h_2(t) \sqrt{\int_{\I}(\chi^2 - g(t))^2}$$
where:
$$h_2(t) = \sqrt{\int_{\I} \abs{v_x(t,x)}^2\dd x}$$
\item 
Notice that since for all $t\in[0,1]$ $\chi(t,\cdot)\in\Sob{1}(\I)$ we have:
\begin{equation*}
\begin{split}
\chi^2(t,x)-\chi^2(t,x_0) &= \int_{x_0}^{x}(\chi^2)_x(t,y)\dd y \\
 &=  2\int \chi(t,y)\chi_x(t,y)\dd y \\ 
 &\leq 2\sqrt{\int_{\I}\chi^2}\sqrt{\int_{\I} (\chi_x)^2}
\end{split}
\end{equation*}
By integrating first with respect to $x_0$, squaring and then integrating with respect to $x$ we get:
$${\int_{\I}(\chi^2(t,x) - g(t))^2}\dd x \leq 
4{gh}$$
so from \eqref{vxchisquare-estimate} we obtain:
\begin{equation}\label{vxchisquare-estimate2}
\int_{\I} v_x \chi^2 \leq h_2(t)\cdot 2\cdot \sqrt{g \cdot h} + h_1(t)\cdot g(t)
\end{equation}
where
$$h_1(t) = \int_{\I} \abs{v_x(t,x)}\dd x$$
\item 
By Young's inequality:
$$2h_2\sqrt{gh}  \leq \frac{g\cdot h_2^2}{2\mu} + 2\mu h$$
\item \eqref{vxchisquare-estimate2} now becomes:
$$\frac 12 \int_{\I} v_x \chi^2 \leq \frac{g\cdot h_2^2}{4\mu} + \mu h(t) + \frac 12 h_1(t)\cdot g(t)$$
\end{enumerate}
\item Combining the last estimate with \eqref{pairingphiminus-gh} yields:
$$g'(t) \leq  h_3(t)\cdot g(t)$$
where
$$h_3(t) = \left(\frac{h_2^2(t)}{2\mu} + h_1(t)\right)$$
\item $h$ is integrable and $g(0) = 0$ so $g(t) = 0$ for any $t>0$.
\end{enumerate}
\end{pr}

\begin{prop}\label{Kequalszero-lma}
If $([\varphi],K)\in\Solthree$  then $K=0$.
\end{prop}
\begin{pr}
 From any representative $\varphi$  we define $\psi\in\Sob{1,2}(\IxI)$ by $\psi(t,x) = \e^{-Kt}\varphi(t,x)$. A simple computation shows that $\psi_t - \mu\psi_{xx} + v\psi_x = 0$ on $\IxI$ so we may use the preceding lemma to get that the trace $\psi(0)\in\Cont{0}(\overline{\I})$ is well defined. We denote its minimum and maximum values in the following way:
\begin{align*}
\gamma_+ &= \max_{x\in \I} \psi(0,x) \\
\gamma_- &= \min_{x\in \I} \psi(0,x) 
\end{align*}
so
\begin{align*}
\psi_1:=\gamma_+ - \psi &\geq 0\\
\psi_2:=\psi-\gamma_- &\geq 0
\end{align*}
$\psi_1$ and $\psi_2$ both qualify for the preceding lemma and $\psi_1(0)\geq 0$ and $\psi_2(0)\geq 0$. Moreover by construction we have $\psi(1) = \e^{-K}\psi(0)$ so applying the lemma yields:
\begin{align*}
\psi_1(1) = \gamma_+ - \e^{-K} \psi(0) &\geq 0\\
\psi_2(1) = \e^{-K} \psi(0) - \gamma_- &\geq 0
\end{align*}
so we get $\e^{-K}\leq 1$ and $\e^{-K}\geq 1$ and hence $K=0$.
\end{pr}

 So the ``eigenvalue'' $K$ must be zero. We now prove that  degeneracy for the corresponding eigenspace is impossible, i.e. we prove that it is one dimensional. Degeneracy implies indeed that the eigenspace would meet the boundary of the cone of positive functions. We will show that this cannot occur because 
 $\Soltwo$ is bounded in $\Cont{0}(\overline{\Q})/\R$. This fact, which follows from our a priori estimate of \autoref{aprioriestimatethm}, is proved in the next lemma.
 
 Notice that by the previous propositions $\Soltwo$ is naturally embedded in $\HWq$ (instead of $(\HWq)\times\R$) which is itself embedded in $\Cont{0}(\overline{\Q})$. With this identification we have the following lemma:

\newcommand{\ImL}{L}

\begin{lma}
$\Soltwo$ is bounded in $\Cont{0}(\overline{\Q})/\R$.
\end{lma}
\begin{pr}
\begin{enumerate}
\item We first show that $\Linop$ (defined by $\Linop u = u_t -\mu u_{xx}$ for $u\in\HWq$) is an isometry from $\HWq$ to $\ImL :=\bigl\{g\in\Leb{2}(\Q);\ \int_\Q g = 0\bigr\}$. Take $f\in\ImL$. We define $\tilde{f}$ on $\T\times(-1,1)$ by symmetrisation:
$$\begin{cases}\tilde{f}(t,x) = f(t,x)&\text{if $x \geq 0$}\\
\tilde{f}(t,x) = f(t,-x) & \text{if $x\leq 0$}
\end{cases}
$$
\newcommand{\doubletorus}{\T\times(\R/2\Z)}
We can now regard $\tilde{f}$ as an element of $\Leb{2}(\doubletorus)$. Notice that $\int_{\doubletorus} \tilde{f} = 0$. By Fourier analysis there is a unique $\tilde{u}\in\Sob{1,2}(\doubletorus)/\R$ solution of 
$$\tilde{u}_t - \mu\tilde{u}_{xx} = \tilde{f}$$
Now $\tilde{u}_x\in\Sob{\frac 12, 1}(\doubletorus)$ so $u_x$ has a trace on $\T\times\{0\}$ and $\T\times\{1\}$. By symmetry it must be zero in $\Leb{2}(\T)$.
\item 
For $(W,0)\in\Soltwo$ we have: 
$$\Linop W = -\frac 12 (W_x)^2 - v W_x $$
so by the previous result: 
$$W = -\frac{1}{2}\Linop^{-1}\bigl(W_x(2v+W_x)\bigr)$$
\item Since $W\in \Soltwo \implies W_x \in  \Solone$ (cf. \autoref{cardinalityprop}), and $W_x\in\Solone \iff \Burgop(v+W_x) = \Burgop(v)$, by the a priori estimate (\autoref{aprioriestimatethm})  there exists $C>0$ such that $\forall W\in\Soltwo$, $\norm{W_x}\leq C$. Using $\H\subset\Leb{4}$ we obtain that $W_x(2v+W_x)$ is bounded in $\Leb{2}$. 
\item Combining the three preceding steps we conclude that $\Soltwo$ is bounded in $\HWq$. But $\HWq\subset\Cont{0}(\overline{\Q})/\R$ so we get the result.
\end{enumerate}
\end{pr}

\begin{prop}\label{descriptionSolthreeLma}
$$\Solthree = \bigl\{([1],0)\bigr\}$$
\end{prop}
\begin{pr}
It is obvious that $([1],0)\in\Solthree$. 
We proved (in \autoref{Kequalszero-lma}) that $(\varphi,K)\in\Solthree$ implies $K=0$. So we take $([\varphi],0)\in\Solthree$ and we will now show that $[\varphi]=[1]$.
\begin{enumerate}
\item
Given one representative $\varphi$ of $[\varphi]$, let us choose a point $x_0\in\overline{\Q}$ where $\varphi$ takes its minimum $\gamma$: 
$$\gamma = \min_{\overline{\Q}} \varphi = \varphi(x_0)$$
 and let us define for $n\in\N$ the function $\psi_n$ on $\overline{\Q}$ by
 $$\psi_n := {\varphi} - \gamma+\frac{1}{n}$$
By construction we have $ \psi_n(x_0) \xrightarrow[n\rightarrow\infty ]{}0$ and thus $\log \psi_n (x_0)\xrightarrow[n\rightarrow\infty ]{}\infty $. It is also clear that for any $n\in\N$ we have $\psi_n\in\Solthree$ so by \autoref{exploglma}, $[\log \psi_n]\in\Soltwo$.  
\item If we now assume that $[\varphi]\neq[1]$ then there exists $x_1$ such that ${\varphi}(x_1)\neq\gamma$ so the sequence $\{\log(\psi_n(x_1))\}_{n\in\N}$ is bounded. As a result $\norm{\log(\psi_n)}_{\Cont{0}/\R}\xrightarrow[n\rightarrow\infty ]{}\infty $. That is a contradiction to the fact that $\Soltwo$ is bounded in $\Cont{0}(\overline{\Q})/\R$. We may therefore conclude that $[\varphi] = [1]$ and the proposition is proved.
\end{enumerate}
\end{pr}

\section{Smoothness}

At this point we have all the ingredients to prove our main result, the \autoref{SmoothnessThm}:

\begin{pr}[Proof of \autoref{SmoothnessThm}]
\begin{enumerate}
\item $\Burgop'(m)$ is injective for all $m\in\H$ because of the observation \eqref{square-derivative-observation} and \autoref{uniquenessThm}.
\item $\Burgop'(m) = \Linop + S'(m)$ but $S'(m)$ is compact (\autoref{Scompact}) so the Fredholm index $\ind\bigl(\Burgop'(m)\bigr) = \ind(\Linop) = 0$ since $\Linop$ is an isomorphism (cf. \autoref{linearcasethm}). Thus $\Burgop'(m)$ is surjective. Since $\Burgop'(m)$ is continuous, linear and bijective it is a homeomorphism.
\item We can use the inverse mapping theorem in Banach spaces (see \cite{Hormander}) to assert that $\Burgop$ is locally a diffeomorphism.
\item $\Burgop$ is a surjection (\autoref{nonlinearsurjectivity}) so it is a global diffeomorphism.
\end{enumerate}
\end{pr}

\section{Extensions}
Our method can be adapted to cover the case of non homogeneous Dirichlet boundary conditions as well as the case of periodic spatial boundary conditions (prescribing the momentum  $\int u(t,x)\dd x$ at $t=0$ for the solution).
%



\begin{thebibliography}{12}
\providecommand{\natexlab}[1]{#1}
\providecommand{\url}[1]{\texttt{#1}}
\expandafter\ifx\csname urlstyle\endcsname\relax
  \providecommand{\doi}[1]{doi: #1}\else
  \providecommand{\doi}{doi: \begingroup \urlstyle{rm}\Url}\fi


\bibitem{Brezis}
Ha{\"{\i}}m Brezis.
\newblock \emph{Analyse fonctionnelle}.
\newblock Collection Math\'ematiques Appliqu\'ees pour la Ma\^\i trise.
  [Collection of Applied Mathematics for the Master's Degree]. Masson, Paris,
  1983.
\newblock ISBN 2-225-77198-7.
\newblock Th\'eorie et applications. [Theory and applications].

\bibitem{Burgers}
J.~M. Burgers.
\newblock Correlation problems in a one-dimensional model of turbulence. {I}.
\newblock \emph{Nederl. Akad. Wetensch., Proc.}, 53:\penalty0 247--260, 1950.

\bibitem{Cole}
Julian~D. Cole.
\newblock On a quasi-linear parabolic equation occurring in aerodynamics.
\newblock \emph{Quart. Appl. Math.}, 9:\penalty0 225--236, 1951.
\newblock ISSN 0033-569X.

\bibitem{Deimling}
Klaus Deimling.
\newblock \emph{Nonlinear functional analysis}.
\newblock Springer-Verlag, Berlin, 1985.
\newblock ISBN 3-540-13928-1.

\bibitem{WeinanE}
Weinan E.
\newblock Aubry-{M}ather theory and periodic solutions of the forced {B}urgers
  equation.
\newblock \emph{Comm. Pure Appl. Math.}, 52\penalty0 (7):\penalty0 811--828,
  1999.
\newblock ISSN 0010-3640.

\bibitem{Fokas}
A.~S. Fokas and J.~T. Stuart.
\newblock The time periodic solution of the {B}urgers equation on the half-line
  and an application to steady streaming.
\newblock \emph{J. Nonlinear Math. Phys.}, 12\penalty0 (suppl. 1):\penalty0
  302--314, 2005.
\newblock ISSN 1402-9251.

\bibitem{Fontes-Saksman}
Magnus Fontes and Eero Saksman.
\newblock Optimal results for the two dimensional {N}avier-{S}tokes equations
  with lower regularity on the data.
\newblock In \emph{Actes des Journ\'ees Math\'ematiques \`a la M\'emoire de
  Jean Leray}, volume~9 of \emph{S\'emin. Congr.}, pages 143--154. Soc. Math.
  France, Paris, 2004.

\bibitem{Hopf}
Eberhard Hopf.
\newblock The partial differential equation {$u\sb t+uu\sb x=\mu u\sb {xx}$}.
\newblock \emph{Comm. Pure Appl. Math.}, 3:\penalty0 201--230, 1950.
\newblock ISSN 0010-3640.

\bibitem{Hormander}
Lars H{\"o}rmander.
\newblock \emph{The analysis of linear partial differential operators. {I}}.
\newblock Classics in Mathematics. Springer-Verlag, Berlin, 2003.
\newblock ISBN 3-540-00662-1.
\newblock Distribution theory and Fourier analysis, Reprint of the second
  (1990) edition [Springer, Berlin; MR1065993 (91m:35001a)].

\bibitem{Moser}
H.~R. Jauslin, H.~O. Kreiss, and J.~Moser.
\newblock On the forced {B}urgers equation with periodic boundary conditions.
\newblock In \emph{Differential equations: La Pietra 1996 (Florence)},
  volume~65 of \emph{Proc. Sympos. Pure Math.}, pages 133--153. Amer. Math.
  Soc., Providence, RI, 1999.

\bibitem{Kreiss}
Heinz-Otto Kreiss and Jens Lorenz.
\newblock \emph{Initial-boundary value problems and the {N}avier-{S}tokes
  equations}, volume 136 of \emph{Pure and Applied Mathematics}.
\newblock Academic Press Inc., Boston, MA, 1989.
\newblock ISBN 0-12-426125-6.

\bibitem{Sinai}
Ya.~G. Sina{\u\i}.
\newblock Two results concerning asymptotic behavior of solutions of the
  {B}urgers equation with force.
\newblock \emph{J. Statist. Phys.}, 64\penalty0 (1-2):\penalty0 1--12, 1991.
\newblock ISSN 0022-4715.

\end{thebibliography}
\end{document}